\documentclass[11pt,reqno]{amsart}      
\usepackage{amssymb,amsthm}             
\numberwithin{equation}{section}        

\renewcommand{\Re}{{\mathbb R}}         
\newcommand{\half}{\frac{1}{2}}         
\newcommand{\third}{\frac{1}{3}}        
\newcommand{\Eps}{\epsilon}             
\newcommand{\starW}{{^*W}}              
\newcommand{\BR}{Q}                     
\newcommand{\EBR}{{\mathcal Q}}         
\newcommand{\FF}{\mathcal F}            
\newcommand{\Lie}{\mathcal L}		
\newcommand{\tr}{\text{\rm tr}}		

\newcommand{\Ric}{\text{Ric}}

\renewcommand{\div}{\text{div}}
\newcommand{\Lapse}{N}
\newcommand{\Shift}{X}
\newcommand{\curl}{\text{\rm curl}}
\newcommand{\aM}{\bar M}              
\newcommand{\anabla}{\bar \nabla}     
\newcommand{\ame}{\bar g}             
\newcommand{\agamma}{\bar \gamma}     
\newcommand{\aR}{\bar R}              

\theoremstyle{plain}
\newtheorem{thm}{Theorem}[section]

\newtheorem{remark}{Remark}[section]

\title[Bel-Robinson energy]{Bel--Robinson energy and constant mean curvature
foliations}
\author[L. Andersson]{Lars Andersson}
\thanks{Supported in part by the Swedish Research Council, 
contract no.  R-RA 4873-307, the NSF,
contract no. DMS 0104402, and the Erwin Schr\"odinger Institute, Vienna.}
\address{Department of Mathematics\\
University of Miami\\
Coral Gables, FL 33124\\
USA}
\email{larsa\char'100math.miami.edu}

\date{July 20, 2003}

\allowdisplaybreaks[2]

\begin{document}
\begin{abstract} An energy estimate is proved 
for the Bel--Robinson energy along a
constant mean curvature foliation in a spatially compact vacuum spacetime, 
assuming an $L^{\infty}$ bound on the second fundamental form, and a bound on
a spacetime version of Bel--Robinson energy. 
\end{abstract} 
\maketitle
\section{Introduction}
Let $(\aM, \ame)$ be a 3+1 dimensional $C^{\infty}$ 
maximal  globally hyperbolic vacuum (MGHV) 
space--time, which is spatially compact, i.e. $\aM$ has compact Cauchy
surfaces. 
One of the main conjectures (the CMC conjecture, see \cite{andersson:survey}
for background) concerning spatially compact MGHV spacetimes
states that if there is a constant mean curvature (CMC) Cauchy surface $M_0$, 
in such a spacetime $\aM$, then there is a foliation in $\aM$ of CMC Cauchy
surfaces with mean curvatures taking on all geometrically 
allowed values. 
Specifically, in case the Cauchy surface $M_0$ is of Yamabe type 
$-1$ or $0$, then the mean curvatures take all values in $(-\infty, 0)$, 
or $(0,\infty)$, depending on the sign of the mean curvature of $M_0$, 
while in case $M_0$ is of Yamabe type $+1$, the mean curvatures take on 
all values in $(-\infty, \infty)$. 
 
The only progress towards proving the CMC conjecture so far has been made
under conditions of symmetry, cf.  \cite{rendall:cosmologicalCMC}, 
or curvature bounds
\cite{anderson:longtime}, \cite{andersson:cmcflat}. 

One approach to the CMC conjecture is to view it as a statement about the
global existence problem for the Einstein vacuum field equations 
\begin{equation}\tag{EFE} 
\aR_{\alpha\beta} = 0 ,
\end{equation} 
in the CMC time gauge. It is known that in the CMC gauge with zero shift, 
the (EFE) are non--strictly hyperbolic \cite{ChB:ruggeri} while in
other gauges such as wave coordinates, the (EFE) form a system of 
quasi--linear wave equations for the metric $\ame$. In this context, it has
been conjectured that the Cauchy problem for the (EFE) is 
well--posed for data in $H^2\times H^1$ (the $H^2$ conjecture, see 
\cite{klainerman:ucla}).

From this point of view it is interesting to consider continuation principles
for the (EFE), in CMC gauge. In this note we will use a scaling 
argument to prove an energy
estimate for CMC foliations. The energy we consider is a version of 
the Bel--Robinson
energy. For a spatial hypersurface $M$ in $\aM$, the energy expression we
consider is defined by 
$$
\EBR(M) = \int_M (|E|^2 + |B|^2 )\mu_g ,
$$
where $E,B$ are the electric and magnetic parts of the Weyl tensor (defined
w.r.t. the timelike normal $T$ of $M$. Roughly speaking, $\EBR$ bounds 
Cauchy data $(g,K)$ on $M$ in $H^2\times H^1$. Here, 
$g$ is the induced metric on $M$ and $K$ is the
second fundamental form of $K$.   Therefore, if the $H^2$ conjecture is true,
apriori bounds for the Bel--Robinson energy can be expected to be relevant
to the global existence problem for the (EFE). 

Let $H = \tr K$ denote the mean curvature, and assume the CMC gauge condition
$H = t$.
Define the spacetime Bel--Robinson energy of a CMC foliation 
$\FF_I = \{M_t, t \in I\}$ by 
\begin{equation}\label{eq:QF}
\EBR(\FF_I) = \int_I dt \int_{M_t} \Lapse (|E|^2 + |B|^2) \mu_g ,
\end{equation}
where $\Lapse$ is the lapse function. 

We are now ready to state our main result
\begin{thm}\label{thm:main} Let $(\aM,
\ame)$ be a  MGHV space--time, and let 
$I = (t_-,t_+)$ with  $-\infty < t_- < t_+ < 0$, be such that
there is a CMC foliation  $\FF_I$ in $(\aM,\ame)$. Let $t_0 \in I$.

Suppose that $\limsup_{t \to t_+} \EBR(t) = \infty$. Then at least one of the
following holds: 
\begin{enumerate}
\item \label{point:EBR}
$\limsup_{t \nearrow t_+} \EBR(\FF_{[t_0, t)})  = \infty$,
\item \label{point:K}
$\limsup_{t \nearrow t_+} \frac{||K(t)||_{L^{\infty}}}{|H(t)|} = \infty$.
\end{enumerate}
The time reversed statement with $t_+$ replaced by $t_-$ also holds. 
\end{thm}

\begin{remark} Let $(M,\gamma)$ be a compact hyperbolic 3--manifold with
sectional curvature $-1$. Then the metric $\agamma = - d\rho^2 + \rho^2
\gamma$ on $\aM = (0,\infty) \times M$ is flat. It follows from the work of
Andersson and Moncrief \cite{andersson:moncrief:global} that for small
perturbations of $(\aM, \agamma)$, there is a global CMC foliation $\FF_{[t_0,t)}$ in the
expanding direction, and for this foliation, the Bel--Robinson energy decays
as $\EBR(t) = O (H^2(t))$, which implies that the space--time Bel--Robinson
energy $\EBR(\FF)$ is bounded in this case. It is interesting to consider the
behavior of $\EBR(\FF_{[t_0,t)})$ when $t_0$ decreases. 
\end{remark}

The proof of Theorem \ref{thm:main} is based on a scaling argument, which we
now sketch. The statement of the theorem is symmetric in time, but here we
consider only the future time direction, the argument in the reverse
direction is similar. 
Suppose for a contradiction there is a constant $\Lambda <
\infty$ so that 
$ \EBR(\FF_{[t_0, t_*)}) \leq \Lambda$, 
$||K||_{L^{\infty}}^2/H^2 \leq \Lambda$ for all $t \in [t_0,t_*)$, and 
that 
$$
\limsup_{t \nearrow t_*} \EBR(t) = \infty \, .
$$
An energy estimate shows that $r_h(t)\EBR(t) \leq C$, where $r_h$ is an
$L^{1,p}$ harmonic radius, for some fixed $p$, $3 < p < 6$, 
and hence if $r_h$ is bounded from
below there is nothing to prove. Suppose for a contradiction that $r_h \to 0$
as $t \nearrow t_*$. 
The combination $r_h \EBR$ is scale invariant, and hence by
rescaling $\ame$ to $\ame' = r_h^{-2} \ame$, 
we get a sequence of metrics $g'$ with $\EBR$
bounded. $\EBR$ bounds $g'$ in $L^{2,2}$ and hence we may 
pick out a subsequence of
$(g',\Lapse')$, 
which converges weakly to a solution $(g_{\infty},\Lapse_{\infty})$ 
of the static vacuum Einstein equations,
\begin{subequations}\label{eq:static}
\begin{align}
\Delta \Lapse &= 0 , \\
\nabla^2 \Lapse &= \Lapse \Ric .
\end{align}
\end{subequations}
It follows from our assumptions that the limit $g_{\infty}$ is complete, 
and the
limiting $\Lapse_{\infty}$ is bounded from 
above and below. Then by \cite{anderson:staticI}, 
$g_{\infty}$ must be flat, with infinite harmonic radius, 
which contradicts $r'_h = 1$, by the weak continuity of $r_h$ on $L^{2,2}$.
We conclude that in fact $r_h$ is bounded away from zero. and hence
that $\EBR$ does not blow up, which 
proves the
theorem.

\section{Preliminaries}
For a space--like hypersurface $M$ in $\aM$ we denote its 
timelike normal $T$ and induced metric and second fundamental form
$(g,K)$. We assume all fields are $C^\infty$ unless otherwise stated.  
Let lower case greek indices run over $0,\dots,3$ while lower case latin
indices run over $1,\dots,3$. 
We work in an adapted frame $e_{\alpha}$, with $e_0 = \partial_t$. 
Our convention for $K$ is $K_{ab} = -\half \Lie_T \ame_{ab}$, 
so that if the mean curvature $H=\tr K$ is negative,
$T$ points in the expanding direction. We will sometimes 
use an index $T$ to denote contraction with
$T$, for example $u_{T} = u_{\alpha}T^{\alpha}$. 

In a nonflat spatially compact, globally hyperbolic, vacuum
spacetime, the maximum principle implies 
uniqueness of constant mean curvature (CMC) 
Cauchy surfaces. In particular,  each $x \in \aM$ is
contained in at most one CMC Cauchy surface, and for each $t \in \Re$, there is
at most one $M_t$ with mean curvature $t$. 

Let $I \subset \Re$ be an interval. 
A CMC foliation $\FF_I$ in $\aM$ is a  foliation $\FF_I = \{ M_t , t \in
I\}$ such that for each $t \in I$, $M_t$ is a $C^{\infty}$ 
CMC Cauchy surface with mean curvature $t$. When convenient we will write
$g(t), K(t)$ for the data induced on $M_t$. Introducing coordinates
$x^{\alpha}$ with $x^0 = t$, the lapse and shift $\Lapse,\Shift$ of the
foliation are defined by $\partial_t = \Lapse T + \Shift$. We may without
loss of generality assume $\Shift = 0$.  

We call $\FF_I$ a maximal CMC
foliation in $\aM$ if there is no interval $I'$ containing $I$ as a strict
subset with a CMC foliation $\FF_{I'}$. Given a foliation in $\aM$, we write
$\aM_{\FF}$ for the support of $\FF$.

Assume that $\aM$ contains a compact, 
constant mean curvature (CMC) Cauchy surface $M_0$ with  mean curvature 
$H^0 < 0$. 
By standard results 
there is then 
an interval $I = (t_-, t_+) \subset \Re$, 
$H^0 \in I$, such that there is a CMC
foliation $\FF_I$, and by uniqueness, $M_{H^0} = M_0$. 
Hence if $\FF_I$ is a maximal CMC
foliation, then $I$ is open. 

\subsection{The Bel--Robinson energy}
Let $W$ be the Weyl tensor of $(\aM, \ame)$ and let $\starW$ denote its
(left) dual (in vacuum $\starW = W^*$). 
The Bel--Robinson tensor $\BR$ of $(\aM, \ame)$ is defined by 
\begin{equation}
\begin{split}
\BR_{\alpha\beta\gamma\delta} &=  W_{\alpha\mu\gamma\nu}W_{\beta\ \delta}^{\ \mu\ \nu} + 
\starW_{\alpha\mu\gamma\nu}\starW_{\beta\ \delta}^{\ \mu\ \nu} . 
\end{split}
\end{equation}
Then $\BR$ is totally symmetric and trace--less, and in vacuum, $\BR$ has
vanishing divergence. 

Let $E_{\alpha\beta} = W_{\alpha T \beta T}$,
$B_{\alpha\beta} = \starW_{\alpha T \beta T}$
be the electric and magnetic parts of the Weyl tensor. Then $E,B$ are
symmetric, 
$t$--tangent (i.e. 
$E_{\alpha T} = B_{\alpha T} = 0$) and trace invariant,
$g^{ab} E_{ab} = g^{ab} B_{ab} = 0$. 

In vacuum, we have 
\begin{align*}
E_{ab} &= \Ric_{ab} + H K_{ab} - K_{ac}K^c_{\ b} ,\\
B_{ab} &= - \curl K_{ab} ,
\end{align*}
where for a symmetric tensor in dimension $3$, 
$$
\curl A_{ab} = \half ( \Eps_a^{\ st} \nabla_t A_{sb} + \Eps_b^{\ st}
\nabla_t A_{sa}) .
$$
Recall that for symmetric traceless tensors in dimension $3$, the Hodge system 
$A \mapsto (\div A, \curl A)$ is elliptic.  

The following identities, see \cite{andersson:moncrief:global}, 
relate $\BR$ to $E$ and $B$, 
\begin{subequations}\label{eq:BRid}
\begin{align}
\BR_{TTTT} &= E_{ab}E^{ab} + B_{ab} B^{ab} = |E|^2 + |B|^2 ,
 \label{eq:BREB}\\
\BR_{a TTT} &= 
2 (E \wedge B)_a , \label{eq:BREBi000} \\
\BR_{abTT} &= 
-  (E \times E)_{ab} -  ( B \times B)_{ab} 
+ \third ( |E|^2 + |B|^2 ) g_{ij} ,
\label{eq:BREBij00}
\end{align}
\end{subequations}
where by definition, for symmetric tensors $A,B$ in dimension 3, 
\begin{align*}
(A \wedge B)_a &= \Eps_a^{\ bc} A_b^{\ d}B_{dc}, \\  
(A \times B)_{ab} &= \Eps_a^{\ cd} \Eps_b^{\ ef} A_{ce} B_{df} + \third (A
\cdot B) g_{ab} - \third (\tr A)(\tr B) g_{ab} . \\
\end{align*} 
From equation (\ref{eq:BREB}) it follows that 
$\BR_{TTTT}\geq 0$ with equality if and only if $W = 0$. 
Let $\FF$ be a foliation in $\aM$. 
The 
Bel--Robinson energy $\EBR(t)$ 
of $M_t \in \FF$ w.r.t. the time--like normal $T$, 
is defined by
$$
\EBR(t) =\EBR(M_t) 
= \int_{M_t} \BR_{TTTT} \mu_g = \int_{M_t} (|E|^2 + |B|^2)\mu_g .
$$
An application of the Gauss law gives in vacuum, 
$$
\partial_t \EBR(t) = - 3 \int_{M_t}
\Lapse \BR_{\alpha\beta TT}\pi^{\alpha\beta}
\mu_g ,
$$
where $\pi_{\alpha\beta} = \anabla_{\alpha} T_{\beta}$. 
A computation shows that the only nonzero components of $\pi_{\alpha\beta}$ 
are $\pi_{ab} =  - K_{ab}$, $\pi_{T a} = \Lapse^{-1} \nabla_a \Lapse$.
Thus
\begin{equation}\label{eq:lapseBR}
\Lapse \BR_{\alpha\beta TT}\pi^{\alpha\beta} = - \Lapse \BR_{abTT} K^{ab} 
+  \BR_{aTTT} \nabla^a \Lapse .
\end{equation}

\section{Proof of Theorem \protect{\ref{thm:main}}}
We will assume that the complement of points \ref{point:EBR},\ref{point:K} of
Theorem \ref{thm:main} holds, and prove from this that $\EBR(t)$ does not
blow up.  
Assume for a contradiction there is a 
constant $\Lambda > 1$ so that for $t \in [t_0,t_*)$, 
\begin{equation} \label{eq:Lambdabound}
\EBR(\FF_{[t_0, t)}) \leq \Lambda , \qquad 
\frac{||K(t)||_{L^{\infty}}^2}{H^2(t)} \leq \Lambda, 
\end{equation}
and that $\limsup_{t \nearrow t_*} \EBR(t) = \infty$.

We let $L^{s,p}$ denote the $L^p$ Sobolev spaces and write $H^s$ for
$L^{s,2}$. We will
sometimes use subindices $x$ or $t,x$ to distinguish function spaces defined
w.r.t. space or space--time. 
 
For a foliation
$\FF_I$, we may without loss of generality assume that $T = \Lapse^{-1}
\partial_t$, where $\Lapse > 0$ is the lapse function of the foliation. Then
$\ame$ is of the form 
$$
\ame = - \Lapse^2 dt^2 + g_{ab} dx^a dx^b .
$$
The lapse function satisfies 
\begin{equation}\label{eq:lapse}
- \Delta \Lapse + |K|^2 \Lapse = 1 ,
\end{equation}
which using the maximum principle implies the estimate 
\begin{equation}\label{eq:lapse-est}
1/||K||_{L^{\infty}}^2 \leq \Lapse \leq 3/H^2 .
\end{equation}
Let $\FF$ be a foliation in $\aM$. 
From (\ref{eq:lapseBR}) we get 
$$
| \partial_t \EBR | \leq C_1( ||\nabla \Lapse ||_{L^{\infty}} +
|| N ||_{L^{\infty}}||K||_{L^{\infty}}) \EBR ,
$$
with $C_1$ a universal constant.
By assumption, $|K|^2/H^2 \leq \Lambda$. Let $\hat K = K - (H/3) g$ be the
traceless part of $K$. 
Using $|K|^2 = H^2/3 + |\hat K|^2 \geq H^2/3$, we get 
\begin{equation}\label{eq:dtEBR-est-full}
| \partial_t \EBR | \leq C( ||\nabla \Lapse ||_{L^{\infty}} +
\frac{\Lambda}{|H|}) \EBR .
\end{equation}
We may assume without loss of generality that $||\nabla \Lapse
||_{L^{\infty}} \geq \Lambda/|H|$, since otherwise there would be nothing to
prove. Therefore we may absorb $\Lambda/|H|$ in the constant in 
(\ref{eq:dtEBR-est-full}) to get 
\begin{equation}\label{eq:dtEBR-est}
| \partial_t \EBR | \leq C ||\nabla \Lapse ||_{L^{\infty}} \EBR .
\end{equation}

\subsection{The blowup}
Choose once and for all a fixed $p$ satisfying
\begin{equation} \label{eq:12}
3 < p <  6. 
\end{equation}
On the Riemannian manifold $(M_{t}, g_t),$ let $r_{h}(x)$ denote 
the $L^{1,p}$ harmonic radius of $(M_{t}, g_t)$ at $x\in M_{t}$; 
thus $r_{h}(x)$ is the radius of the largest geodesic ball about $x$ on which 
there is a harmonic coordinate chart in which the metric coefficients
$g_{ab}$  satisfy
\begin{equation} \label{eq:13}
r_{h}(x)^{-3/p} ||g_{ab} -  \delta_{ab}||_{L^{p}(B_{x}(r_{h}(x)))}
+ r_{h}(x)^{(p-3)/p} ||\partial g_{ab}||_{L^{p}(B_{x}(r_{h}(x)))}\leq  C, 
\end{equation}
where $C$ is a fixed constant (say $C =$ 1), cf. 
\cite{anderson:rigidity,anderson:cheeger}. By Sobolev 
embedding, in $B_{x}(r_{h}(x)),$ the $C^{\beta}$ norm of $g_{ab}$ is 
controlled, for $\beta  = 1-\frac{3}{p}.$ 
The presence of the 
factors of $r_{h}(x)$ in (\ref{eq:13}) means the estimate (\ref{eq:13}) is 
scale invariant. It follows from this that $r_{h}(x)$ scales as a distance.

It is well known that the Laplacian in such a local harmonic coordinate chart 
on $B_{x}(r_{h}(x))$ has the form
$$
\Delta u = g^{ab}\partial_{a}\partial_{b}u. 
$$
Thus, within $B_{x}(r_{h}(x)), \Delta $ is given in these local coordinates 
as a non--divergence 
form elliptic operator, with uniform $C^{\beta}$ control on 
the coefficients $g^{ab},$ and uniform bounds on the ellipticity constants.

We have the following standard (interior) $L^{p}$ elliptic estimate for this 
Laplace operator, c.f. \cite[Thm. 9.11]{gilbarg:trudinger}. 
Let $B = B_{x}(r_{h}(x))$ and 
$B'  = B_{x}(\frac{1}{2}r_{h}(x)).$ 
Then
\begin{equation} \label{eq:14}
||\Lapse||_{L^{2,p}(B' )} \leq  C(r_{h}(x), p)[||\Delta\Lapse||_{L^{p}(B)} 
+ ||\Lapse||_{L^{p}(B)}]. 
\end{equation}
We drop the dependence on $p$, since $p$ is fixed. We need to make explicit 
the dependence of the constant $C$ on $r_{h}(x).$ This is done by a standard 
scaling argument. Thus, assume (by rescaling if necessary), that 
$r_{h}(x) =1$. Then (\ref{eq:14}) becomes
$$
||\Lapse||_{L^{2,p}(B' )} \leq  C[||\Delta\Lapse||_{L^{p}(B)} +
||\Lapse||_{L^{p}(B)}]. 
$$
By Sobolev embedding, since $p > $ 3 is fixed, and $B'  = B(\frac{1}{2}),$ we 
have
$$
||\nabla\Lapse||_{L^{\infty}(B' )} \leq  c\cdot ||\Lapse||_{L^{2,p}(B' )}, 
$$
so that
$$
||\nabla\Lapse||_{L^{\infty}(B' )}\leq  C_{o}[||\Delta\Lapse||_{L^{p}(B)} 
+ ||\Lapse||_{L^{p}(B)}]. 
$$
and in particular,
\begin{equation} \label{eq:15}
||\nabla\Lapse||_{L^{\infty}(B' )}\leq  
C_{o}[||\Delta\Lapse||_{L^{\infty}(B)} + ||\Lapse||_{L^{\infty}(B)}], 
\end{equation}
where $C_{o}$ is an absolute constant, (i.e. independent of $\Lapse ,$ 
given control on $\Delta $ from definition of $r_{h} =$ 1). Now we put in 
scale factors to make (\ref{eq:15}) scale invariant and write (\ref{eq:15}) as
\begin{equation} \label{eq:16}
r_{h}(x)||\nabla\Lapse||_{L^{\infty}(B' )}\leq  C_{o}[r_{h}(x)^{2}
||\Delta\Lapse||_{L^{\infty}(B)} + ||\Lapse||_{L^{\infty}(B)}]. 
\end{equation}
Note that the function $\Lapse $ is itself scale invariant. Each 
term in (\ref{eq:16}) is invariant under scaling, and thus (\ref{eq:16}) holds in any scale. 
Therefore, it holds in the metric $g(t)$. 

Let 
$$
r_h = r_h(t) = \inf_{s\leq t}\inf_{x \in M_s} r_h(x) .
$$  
From the lapse equation,
$$
\Delta\Lapse = \Lapse|K|^{2} -  1. 
$$
Using (\ref{eq:lapse-est}) and $|K|^2 = H^2/3 + |\hat K|^2$  we find
$$
0 \leq \Delta \Lapse \leq 3 \frac{|\hat K|^2}{H^2} \leq 3\frac{|K|^2}{H^2}
\leq 3 \Lambda .
$$
Thus, we have 
\begin{equation}\label{eq:DN-est}
r_{h}||\nabla\Lapse||_{L^{\infty}(B' )} \leq  
3 C_{o}\left(
r_{h}^{2} \Lambda + \frac{1}{H^2}\right) .
\end{equation}
In particular, this gives the estimate 
\begin{equation}\label{eq:Dlapse-scale}
||\nabla\Lapse(t)||_{L^{\infty}} \leq C(\Lambda,t_*)/r_h(t) . 
\end{equation}

Integrating (\ref{eq:dtEBR-est}), (recall we absorbed the term $\Lambda/|H|$
in (\ref{eq:dtEBR-est-full}) into the constant), gives 
$$
\EBR(t_1) \leq \EBR(t_0) + C \int_{t_0}^{t_1} ds ||\nabla
\Lapse(s)||_{L^{\infty}} \EBR(s) .
$$
We may without loss of generality assume the last term is bigger than 1, so
we may absorb $\EBR(t_0)$ into $C$. 

Multiplying both sides by $r_h$, and
using (\ref{eq:DN-est}) we have
$$
r_h \EBR(t_1) \leq C \int_{t_0}^{t_1} ds \left ( 
r_h^2 \Lambda + \frac{1}{H^2(s)} \right ) \EBR(s) .
$$
We may without loss of generality assume $r_h \leq 1/|H|$, since otherwise there
would be nothing to prove, and therefore we can absorb the term $r_h^2
\Lambda$ into the constant. Then we have 
$$
r_h \EBR(t_1) \leq C \int_{t_0}^{t_1} ds  
\frac{1}{H^2(s)} \EBR(s). 
$$
The inequality (\ref{eq:lapse-est}) implies 
$$
N \geq 3 \Lambda^{-1} H^{-2} , 
$$
which by the definition of the spacetime Bel--Robinson energy 
$\EBR(\FF_{[t_0,t_1)})$, see (\ref{eq:QF}), gives
\begin{equation}\label{eq:integrated-energy-est}
r_h(t_1)\EBR(t_1) \leq C \EBR(\FF_{[t_0,t_1)}) .
\end{equation}
By assumption, $\limsup_{t \nearrow t_*} \EBR(t) = \infty$, which by the assumed
bound on $\EBR(\FF_{[t_0,t_1)})$ 
implies 
$$
\lim_{t \nearrow t_*} r_h(t) =0 .
$$ 
We will show that this contradicts
(\ref{eq:Lambdabound}).

Suppose then that
there is an increasing sequence of  times ${{t_i}}$, 
$t_i \to t_*$ as $i \to \infty$, 
so that $r_i = r_h(t_i)$ satisfy 
$\lim_{i \to \infty} r_i = 0$ (recall that by construction $r_h(t)$ is
decreasing).

Now we have from  (\ref{eq:integrated-energy-est}) and our assumptions,
\begin{equation}\label{eq:BRi}
r_i \EBR(t_i) \leq C .
\end{equation}
Now we introduce the blowup scale. 
Let $\ame'_i = r_i^{-2} \ame$. We will denote the scaled versions of $g,K$ by
$g'_i, K'_i$. We scale 
the coordinates as $t'_i = r_i^{-1} t$, $x'_i= r_i^{-1} x'$, 
so that the coordinate components of 
$g'_i$ are scale invariant.  Then
$|K'_i| \leq \Lambda r_i$, while the lapse $\Lapse$
does not scale, $\Lapse'_i(x'_i) = \Lapse(x)$. 
After translating the time coordinate as $t'_i = r_i^{-1} (t - t_i)$, we 
focus our attention on the time interval $t'_i \in
[-1, 0]$. We further translate the space coordinate so that the center of the
coordinate system (0,0) is the point where the harmonic radius achieves its
minimum value.  

Since $r_h \EBR$ is scale invariant, we have $r_h \EBR =
\EBR'$, and hence the inequality 
\begin{equation}\label{eq:BR-blowup-est}
\EBR'(t'_i) \leq C
\end{equation}
holds. 
This means in view of the definition of the Bel--Robinson energy that at the
blowup scale, $\Ric'_i$ is bounded in $L^2$. 
By construction $r'_h \geq 1$ and  
by \cite{anderson:cheeger}, it follows from the Ricci bound 
that $g'_i$ is bounded in $L^{2,2}_{loc}$. Similarly the Hodge system
relating $K$ to $B$ leads to $K'_i$ bounded in $L^{1,2}_{loc}$. 

The Einstein vacuum equation is scale invariant, and therefore holds at the
blowup scale. We will argue in the next section, 
that the above bounds on $g'_i, K'_i$ allow us to pick out a weakly
convergent subsequence of $(g'_i, \Lapse'_i)$ with limit 
$g_\infty, \Lapse_{\infty}$ solving the static vacuum Einstein equation,
cf. equation (\ref{eq:static}) below, with $g_{\infty}$ complete.

\subsection{Weak convergence}
Let $\ame'_i$ be the sequence of rescaled spacetime metrics. We consider
rescaled time $t'_i$ in the interval $[-1,0]$. By construction, the $L^{1,p}$ 
harmonic radius satisfies $r'_i(0) = 1$, and $r'_i(t) \geq 1$ for 
$t \in [-1,0]$.
Equation (\ref{eq:lapse-est}) implies that the rescaled lapse is bounded from 
above and below, 
\begin{equation}\label{eq:lambda-est-rescaled}
\frac{1}{t_-^2\Lambda} \leq \Lapse'_i \leq \frac{3}{t_+^2} .
\end{equation} 
By (\ref{eq:BR-blowup-est}), we have $\EBR'(t) \leq C$, for $t \in [-1,0]$,
and hence we have $(g'_i,K'_i)$ bounded in $L^{\infty}([-1,0];L^{2,2}_{loc}
\times L^{1,2}_{loc})$. It follows that there is a subsequence which
converges weak-$\star$  to a limit $(g_{\infty}, K_{\infty}) \in 
L^{\infty}([-1,0];L^{2,2}_{loc}\times L^{1,2}_{loc})$, with corresponding
spacetime metric $\ame_{\infty}$. By passing to a further subsequence if
necessary, which we still denote using the index $i$, we may assume that 
$g'_i(0) \rightharpoondown g_{\infty}(0)$ weakly in $L^{2,2}_{loc}$. 

Let us consider the properties of this limit. First note that $|K'_i|\leq
\Lambda r_i \to 0$ as $i \to \infty$, and hence $K_{\infty} \equiv 0$. The
relation $\partial_{t_i'} g'_i = -2\Lapse'_i K'_i$ holds in the limit and
since $K_{\infty} \equiv 0$, we conclude that $g_{\infty}$ is time
independent, so that the limiting spacetime metric $\ame_{\infty}$ is static. 
The lapse equation now implies that the limiting lapse function satisfies 
\begin{subequations}\label{eq:static-Ein-infty}
\begin{equation}\label{eq:lapse-infty}
\Delta_{\infty} \Lapse_{\infty} = 0 ,
\end{equation}
where $\Delta_{\infty}$ is the Laplace operator defined w.r.t. $g_{\infty}$. 

The rescaled spacetime metrics $\ame'_i$ are solutions of the Einstein vacuum 
equation and the evolution equation for $K$, 
$$
\partial_t K = - \nabla^2 N + \Lapse ( \Ric + HK - 2 K : K) ,
$$
where $K : K_{ab} = K_{ac}K_{\ b}^c$, holds weakly in the limit. In view
of the fact that $g \mapsto \Ric$ is weakly continuous on $L^{2,2}_{loc}$ and
$K_{\infty} \equiv 0$, we get the equation 
\begin{equation}\label{eq:Lstatic-evol}
0 = - \nabla_{\infty}^2 \Lapse_{\infty} + \Lapse_{\infty} \Ric_{\infty} .
\end{equation}
\end{subequations}
By construction, $g_{\infty}$ is complete, and hence in view of 
(\ref{eq:static-Ein-infty}) we have a complete solution of the static Einstein 
equations with $N_\infty > 0$. It follows by 
\cite[Theorem 3.2]{anderson:staticI} that $g_{\infty}$ is flat and
$\Lapse_{\infty}$ is constant. In particular, $r_h[g_{\infty}](0) = \infty$. 
Now, since $r_h$ is by definition the
$L^{1,p}$ harmonic radius, $3 < p < 6$, the map $g \to r_h$ is weakly
continuous on $L^{2,2}_{loc}$ and hence by construction 
$r_h[g_{\infty}](0) = 1$. 
This is a contradiction, and it follows that in fact $\liminf_{i \to
\infty} r_i > 0$, which by the BR energy estimate (\ref{eq:BRi}) implies that
$\EBR(t)$ does not blow up. This completes the proof of Theorem
\ref{thm:main}. 
\bigskip

\noindent{\bf Acknowledgements:} 
The problem studied in this paper was suggested by Mike Anderson. I am
grateful to him for many helpful suggestions, and to  
Vince Moncrief and Jim Isenberg for useful discussions
on the topic of this paper.


\providecommand{\bysame}{\leavevmode\hbox to3em{\hrulefill}\thinspace}
\providecommand{\MR}{\relax\ifhmode\unskip\space\fi MR }
\providecommand{\MRhref}[2]{%
  \href{http://www.ams.org/mathscinet-getitem?mr=#1}{#2}
}
\providecommand{\href}[2]{#2}

\end{document}